\documentclass{amsart}

\usepackage{graphicx}         
\usepackage{subfigure} 
\usepackage{url}              

\newtheorem{theorem}{Theorem}[section]

\newtheorem{observation}[theorem]{Observation}

\theoremstyle{definition}
\newtheorem{defn}[theorem]{Definition}

\newcommand{\mfd}{\mathcal{M}}

\newcommand{\regina}{\emph{Regina}}

\newcommand{\snappea}{\emph{SnapPea}}
\newcommand{\snappy}{\emph{SnapPy}}

\newcommand{\Z}{\mathbb{Z}}

\newcommand{\mfda}{\texttt{x101}}
\newcommand{\mfdc}{\texttt{x103}}

\begin{document}

\title{A duplicate pair in the SnapPea census}
\author{Benjamin A.\ Burton}
\address{School of Mathematics and Physics \\
    The University of Queensland \\
    Brisbane QLD 4072 \\
    Australia}
\email{bab@maths.uq.edu.au}
\thanks{The author is supported by the Australian Research Council
    under the Discovery Projects funding scheme (project DP1094516).}

\begin{abstract}
    We identify a duplicate pair in the well-known Callahan-Hildebrand-Weeks
    census of cusped finite-volume hyperbolic 3-manifolds.
    Specifically, the six-tetrahedron non-orientable
    manifolds {\mfda} and {\mfdc} are homeomorphic.
\end{abstract}

\maketitle

%
%

\section{Introduction}

The study of hyperbolic 3-manifolds was transformed by the advent of the
software package {\snappea} in the early 1990s \cite{snappea}.
At the heart of this software is the \emph{cusped hyperbolic census}---a
collection intended to represent all non-compact finite-volume
hyperbolic 3-manifolds that can be
constructed from at most $n$ tetrahedra, for fixed $n$.
Hildebrand and Weeks built the first census in 1989 for
$n=5$ \cite{hildebrand89-cuspedcensusold},
Callahan, Hildebrand and Weeks extended it
to $n=7$ in 1999 \cite{callahan99-cuspedcensus},
and Thistlethwaite recently grew it to $n=8$
in 2010 \cite{thistlethwaite10-cusped8}.

To test whether two cusped hyperbolic 3-manifolds are isometric (and thus
homeomorphic), {\snappea} computes the canonical Epstein-Penner
cell decomposition of each \cite{epstein88-euclidean,weeks93-convex}:
the manifolds are then homeomorphic if and only if their canonical cell
decompositions are combinatorially isomorphic.

There is a problem, however: {\snappea} works with floating point arithmetic,
and is therefore subject to issues such as round-off error and numerical
instability.
These issues impact upon the isometry\,/\,homeomorphism test as follows:
\begin{itemize}
\item If {\snappea} claims that two manifolds \emph{are} homeomorphic then this
claim is reliable.  {\snappea} computes canonical cell decompositions
using local Pachner moves (bistellar flips) \cite{weeks93-convex}, and
so obtaining isomorphic cell decompositions means that {\snappea}
has essentially found a sequence of Pachner moves that transform one manifold
into the other.
\item If {\snappea} claims that two manifolds are \emph{not} homeomorphic,
this claim may be false.  Although {\snappea}
computes cell decompositions that are homeomorphic to
their respective input manifolds, these might
not be the sought-after Epstein-Penner decompositions, and so
two different cell decompositions need not indicate that the two input manifolds
are distinct.
\end{itemize}

As a result, it is theoretically possible that two manifolds with
different names in the cusped hyperbolic census could in fact be homeomorphic.
Here we show that indeed this happens:
the non-orientable manifolds {\mfda} and
{\mfdc}---both described using six tetrahedra in the 1999 census
of Callahan, Hildebrand and Weeks---are in fact the same 3-manifold.

The proof that {\mfda} and {\mfdc} are homeomorphic is
simple: their triangulations from the census are related by just two
Pachner moves (or bistellar flips).  What is much more remarkable
is that this duplicate pair has gone undetected for so long.

In this brief paper we describe the triangulations of
{\mfda} and {\mfdc}, and show why they are homeomorphic.
We also describe the canonical cell decompositions that
{\snappea} computes (at least one of which is erroneous),
in order to to better understand how {\snappea} reaches its
conclusion that the manifolds are distinct.

All computations were performed using the software packages
{\regina} \cite{burton04-regina,regina} and {\snappea} \cite{snappea}
(the latter refers to the mathematical kernel at the heart of both the
original {\snappea} application and its modern Python-based successor
{\snappy} \cite{snappy}).

\section{The duplicate pair}

As is standard nowadays for hyperbolic 3-manifolds, we work with ideal
triangulations.  An \emph{ideal triangulation} of a non-compact
hyperbolic 3-manifold $\mfd$ is essentially a collection of $n$
tetrahedra whose faces are affinely identified in pairs using the
quotient topology, and where
(i)~the link of every vertex is a torus or Klein bottle; and
(ii)~the manifold $\mfd$ is recovered by removing the vertices of each
tetrahedron.
For a richer and more precise discussion on ideal triangulations,
we refer the reader to \cite{thurston86-deformation}.

From here onward we use the labels {\mfda} and {\mfdc} to denote explicit
\emph{triangulations} (since the underlying manifolds are the same).
These triangulations are shipped in the cusped hyperbolic census
with both {\regina} and {\snappy},
and readers can use {\regina} to study their combinatorial structures
in detail.  Here we give a high-level outline of how each
triangulation is constructed.

\begin{figure}[tb]
    \centering
    \subfigure[A five-tetrahedron cube]{%
        \label{fig-cube} \qquad\includegraphics[scale=1.00]{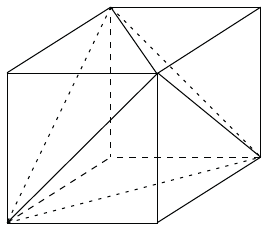}\qquad}
    \hspace{1.5cm}
    \subfigure[Attaching a sixth tetrahedron]{%
        \label{fig-layer} \includegraphics[scale=1.00]{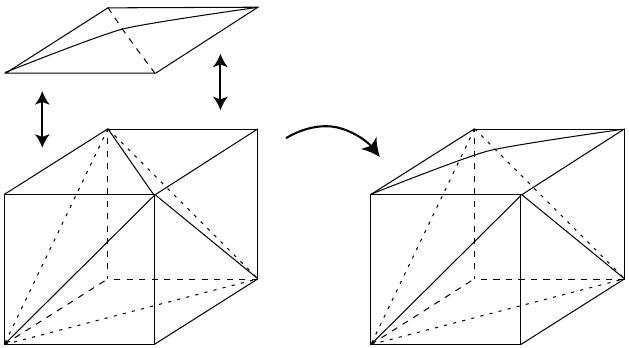}} \\
    \subfigure[Identifying opposite squares]{%
        \label{fig-identify}
        \qquad\includegraphics[scale=1.00]{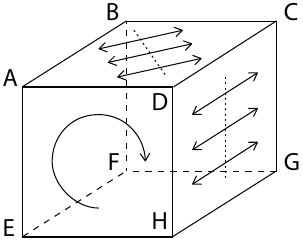}\qquad}
    \caption{Building {\mfda}}
    \label{fig-x101}
\end{figure}

\begin{defn}[Triangulation {\mfda}]
    The six-tetrahedron ideal triangulation {\mfda} is formed as follows.
    \begin{enumerate}
        \item Build a standard five-tetrahedron cube, as shown in
        Figure~\ref{fig-cube}.
        \item Attach an extra tetrahedron to the upper face of the
        cube to change the diagonal on the boundary, as illustrated in
        Figure~\ref{fig-layer}.
        \item Identify opposite faces of the cube, as illustrated in
        Figure~\ref{fig-identify}:
        \begin{itemize}
        \item join the upper and lower squares with a diagonal reflection;
        \item join the left and right squares with a horizontal reflection;
        \item join the front and back squares with a quarter-turn rotation.
        \end{itemize}
        In other words, we
        identify opposite faces according to the following maps:
        \[
        \mathit{ABCD} \longleftrightarrow \mathit{GFEH}\,;\quad
        \mathit{ABFE} \longleftrightarrow \mathit{CDHG}\,;\quad
        \mathit{ADHE} \longleftrightarrow \mathit{CGFB}. \]
    \end{enumerate}
\end{defn}

\begin{figure}[tb]
    \centering
    \subfigure[The degree-four edge $e$]{%
        \label{fig-diagonal} \quad\includegraphics[scale=1.00]{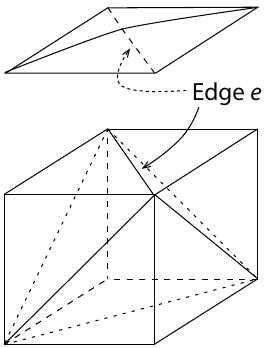}\quad}
    \hspace{0.5cm}
    \subfigure[The octahedron around $e$]{%
        \label{fig-octahedron}
        \quad\ \includegraphics[scale=1.00]{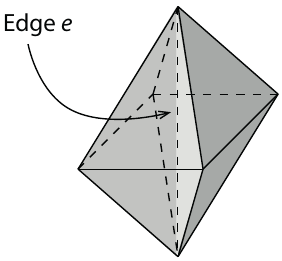}\ \quad}
    \hspace{0.5cm}
    \subfigure[The 4-4 move]{%
        \label{fig-4-4}
        \qquad\includegraphics[scale=1.00]{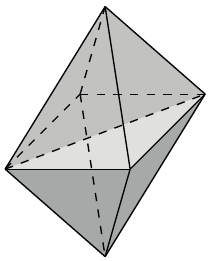}\qquad}
        \\
    \subfigure[The retriangulated cube]{%
        \label{fig-newcube}
        \qquad\includegraphics[scale=1.00]{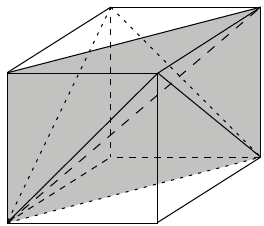}\qquad}
    \caption{Building {\mfdc}}
    \label{fig-x103}
\end{figure}

\begin{defn}[Triangulation {\mfdc}]
    The six-tetrahedron ideal triangulation {\mfdc} is formed as follows.
    Begin with {\mfda}, as described above.  Let $e$ denote the diagonal edge
    between the cube and the extra tetrahedron, as marked in
    Figure~\ref{fig-diagonal}.

    This edge has degree four, meeting four distinct tetrahedra to
    form an octahedron as illustrated in Figure~\ref{fig-octahedron}.
    Replace this octahedron with a different octahedron that has the
    same boundary but a different internal diagonal edge, as shown in
    Figure~\ref{fig-4-4} (this operation is known as a \emph{4-4 move}).

    The final triangulation of the cube is shown in Figure~\ref{fig-newcube};
    note that the shaded rectangle divides the cube into a pair of
    three-tetrahedron triangular prisms.
\end{defn}

Because the identifications between opposite faces of the cube remain
unchanged, the following is immediate:

\begin{observation}
    Triangulations {\mfda} and {\mfdc} are homeomorphic as 3-mani\-folds.
\end{observation}

Moreover, we observe that the two triangulations are ``combinatorially close''.
We can express their relationship in terms of the well-known Pachner moves,
or bistellar flips \cite{pachner91-moves}: since a single 4-4 move can be
expressed as a 2-3 Pachner move followed by a 3-2 Pachner move,
we have:

\begin{observation}
    Triangulations {\mfda} and {\mfdc} are related by just two Pachner moves.
\end{observation}

The common manifold that they represent is non-orientable with first homology
$H_1 = \Z \oplus \Z_2 \oplus \Z_2$.
{\snappea} approximates the hyperbolic volume as $\simeq 5.07470803$.

\begin{figure}[tb]
    \centering
    \subfigure[For \mfda]{%
        \label{fig-ep-layer}
        \includegraphics[scale=1.00]{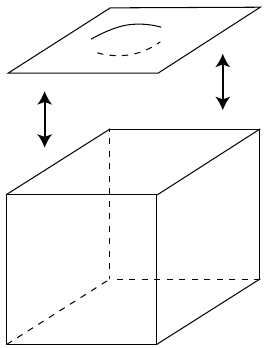}}
    \hspace{1.5cm}
    \subfigure[For \mfdc]{%
        \label{fig-ep-cube}
        \includegraphics[scale=1.00]{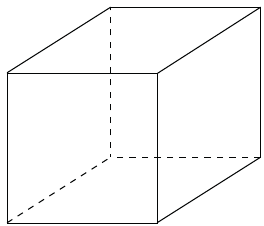}}
    \caption{The ``canonical'' dell decompositions computed by {\snappea}}
\end{figure}

We finish with a brief explanation as to why {\snappea} claims that the
underlying manifolds for {\mfda} and {\mfdc} are distinct.
This is because {\snappea} computes different ``canonical'' cell decompositions
for each:
\begin{itemize}
    \item For {\mfdc}, {\snappea} computes the canonical cell
    decomposition to be a single cube with opposite faces identified, as
    illustrated in Figure~\ref{fig-ep-cube}.

    \item For {\mfda}, {\snappea} computes a decomposition with
    \emph{two cells}: a cube plus a rectangular pillow, joined together
    as illustrated in Figure~\ref{fig-ep-layer}.
\end{itemize}

The cube for {\mfdc} corresponds precisely to the cube illustrated
earlier in Figure~\ref{fig-identify} (i.e., its opposite faces are
identified in the same way).  The cube and pillow for {\mfda} correspond
to the original five-tetrahedron cube plus the sixth ``flat'' tetrahedron
from Figure~\ref{fig-layer}.

The true Epstein-Penner decomposition is indeed the cube found from
{\mfdc} (Weeks gives a simple argument to this effect based on the
symmetry of the cube and its face identifications from
Figure~\ref{fig-identify}).
For {\mfda}, it seems a reasonable conclusion that
the numerical errors that {\snappea} experiences when
computing the canonical cell decomposition are due to the
flat tetrahedron that we attach in Figure~\ref{fig-layer}.

%
%

\bibliographystyle{amsplain}
\bibliography{pure}

\newcommand{\noopsort}[1]{}
\providecommand{\bysame}{\leavevmode\hbox to3em{\hrulefill}\thinspace}
\providecommand{\MR}{\relax\ifhmode\unskip\space\fi MR }
\providecommand{\MRhref}[2]{%
  \href{http://www.ams.org/mathscinet-getitem?mr=#1}{#2}
}
\providecommand{\href}[2]{#2}
\begin{thebibliography}{10}

\bibitem{burton04-regina}
Benjamin~A. Burton, \emph{Introducing {R}egina, the 3-manifold topology
  software}, Experiment. Math. \textbf{13} (2004), no.~3, 267--272.

\bibitem{regina}
Benjamin~A. Burton, Ryan Budney, William Pettersson, et~al., \emph{Regina:
  Software for 3-manifold topology and normal surface theory},
  \texttt{http://\allowbreak regina.\allowbreak sourceforge.\allowbreak net/},
  1999--2013.

\bibitem{callahan99-cuspedcensus}
Patrick~J. Callahan, Martin~V. Hildebrand, and Jeffrey~R. Weeks, \emph{A census
  of cusped hyperbolic 3-manifolds}, Math. Comp. \textbf{68} (1999), no.~225,
  321--332.

\bibitem{snappy}
Marc Culler, Nathan~M. Dunfield, and Jeffrey~R. Weeks, \emph{{SnapPy}, a
  computer program for studying the geometry and topology of 3-manifolds},
  \texttt{http://\allowbreak snappy.\allowbreak computop.\allowbreak org/},
  1991--2013.

\bibitem{epstein88-euclidean}
D.~B.~A. Epstein and R.~C. Penner, \emph{Euclidean decompositions of noncompact
  hyperbolic manifolds}, J. Differential Geom. \textbf{27} (1988), no.~1,
  67--80.

\bibitem{hildebrand89-cuspedcensusold}
Martin~V. Hildebrand and Jeffrey~R. Weeks, \emph{A computer generated census of
  cusped hyperbolic 3-manifolds}, Computers and Mathematics (Cambridge, MA,
  1989), Springer, New York, 1989, pp.~53--59.

\bibitem{pachner91-moves}
Udo Pachner, \emph{P.{L}. homeomorphic manifolds are equivalent by elementary
  shellings}, European J. Combin. \textbf{12} (1991), no.~2, 129--145.

\bibitem{thistlethwaite10-cusped8}
Morwen Thistlethwaite, \emph{Cusped hyperbolic manifolds with 8 tetrahedra},
  \texttt{http://\allowbreak www.\allowbreak math.\allowbreak utk.\allowbreak
  edu/\allowbreak \~{}morwen/\allowbreak 8tet/}, October 2010.

\bibitem{thurston86-deformation}
William~P. Thurston, \emph{Hyperbolic structures on {$3$}-manifolds. {I}.
  {D}eformation of acylindrical manifolds}, Ann. of Math. (2) \textbf{124}
  (1986), no.~2, 203--246.

\bibitem{snappea}
Jeffrey~R. Weeks, \emph{Snap{P}ea: {H}yperbolic 3-manifold software},
  \texttt{http://\allowbreak www.\allowbreak geometrygames.\allowbreak
  org/\allowbreak SnapPea/}, 1991--2000.

\bibitem{weeks93-convex}
\bysame, \emph{Convex hulls and isometries of cusped hyperbolic
  {$3$}-manifolds}, Topology Appl. \textbf{52} (1993), no.~2, 127--149.

\end{thebibliography}

\end{document}